\documentclass{article}%
\usepackage{amsmath}
\usepackage{amsfonts}
\usepackage{amssymb}
\usepackage{graphicx}%
\begin{document}

\title{Some recent trends from research mathematics and their connections to teaching:\\%
Case studies inspired by parallel developments in science and technology}
\author{Palle E. T. Jorgensen\\University of Iowa, Iowa City, IA, USA}
\date{}
\maketitle

\begin{abstract}
\noindent We will outline our ideas for teaching in the core mathematics
disciplines. They are based on our own experience in teaching at a number of
universities in the USA, as well as in Europe. While some of the core ideas
stay and have stayed relatively constant over a long period of time, they must
be varied in accordance with the needs and the demands of students, and they
must constantly updated keeping an eye to current research and to modern
international trends in technology. Our thoughts and suggestions on the use of
these trends in teaching have been tried out by the author, and they are now
in textbooks, some by the author.\medskip


\end{abstract}

\section{\label{Section1}Interaction between technology, research in
mathematics, and teaching}

Mathematics\footnotetext{Expanded version of an invited presentation ``Teaching
of mathematics at various levels in an international
university-system, and connections to research and to current trends in technology''
at the 
Symposium on Mathematics Education Reform
(2006 Beijing, China) http://www.math.ohiou.edu/\symbol{126}shen/calculus/schedule.html .}
draws ideas and strengths from the outside world, and the
connections to parts of engineering have been a boon to mathematics: {}From
signal processing to wavelet analysis!
That is true even if we forget about all of the practical applications
emerging from these connections. Without inspiration from the neighboring
sciences, mathematics would in all likelihood become rather sterile, and
overly formal. I see opportunities at crossroads. Mathematics is reaping
benefits from trends and topics in engineering and in the sciences. It is
witnessed in a striking way by exciting developments in wavelets.

{}From wavelets we see how notions of scale-similarity can be exploited in basis
computations that use tricks devised for signal processing.
This has now all become part of an exciting and fairly recent trend in
mathematics and in technology centered around advances in wavelets and some of
their many applications. But as we outline below, this is only part of a
bigger picture involving also fractal analysis and the use of scale-similarity
in a wide area of scientific problems.

At the same time, the key notion of self-similarity, such as the
scale-similarity used everywhere for wavelets,
is essential to our understanding of fractals: Fern-like pictures that look
the same at small and at large scales.

One problem in the generation of wavelet bases
is selecting the \textquotedblleft nice\textquotedblright\ (here this means
differentiable) wavelets
among huge families of fractal-looking (non-smooth, or singular) functions.

Our analysis will take place in a variety of data sets, or in function spaces.
Here, for the moment we begin with spaces of square-integrable functions
(denoted $L^{2}$); they are Hilbert spaces, and they have special significance
in modeling states in physics, and via the inner product, correlations in
statistics. But $L^{2}$-functions can be very \textquotedblleft
bad\textquotedblright\ indeed!! Computers generate the good and the bad, and
we are left with the task of sorting them out and making selections. It may be
observed (directly from large libraries of pictures \cite{BrJo02b,Jor06a})
that mathematical wavelet
machines are more likely to spit out bad functions unless they are told where
to concentrate the search from the intrinsic mathematics. These wavelets,
signals, and fractals are things that have caught our attention in recent
decades, but the mathematical part of this has roots back at least a hundred
years, for example, to Alfred Haar and to Oliver Heaviside
at the turn of the last century. {}From Haar
we have the first wavelet basis, and with Heaviside
we see the beginning of signal analysis.
It is unlikely that either one knew about the other. Ironically, at the time
(1909), Haar's
paper had little impact and was hardly noticed, even on the small scale of
\textquotedblleft notice\textquotedblright\ that is usually applied to
mathematics papers.

Haar's wonderful wavelet
only began to draw attention in the mid-nineteen-eighties when the connections
to modern signal processing
became much better understood. These connections certainly served as a main
catalyst in what are now known as wavelet
tools in pure and applied mathematics. But at the outset, the pioneers in
wavelets
had to \textquotedblleft rediscover\textquotedblright\ a lot of stuff from
signal processing: frequency bands, high-pass, low-pass, analysis
and synthesis using down-sampling, and up-sampling, reconstruction of signals,
resolution
of images; all tools that have wonderful graphics representations in the
engineering literature. But still, why would we think that Fourier's basis,
and his lovely integral decomposition, are not good enough? Many reasons:
Fourier's method has computational drawbacks. This was less evident before
computers became common and began to play important roles in applied and
theoretical work.

Expansion of functions or signals into basis
decompositions (called \textquotedblleft analysis\textquotedblright
in signal processing) involves basis coefficients (Fourier coefficients,
and so on), and if we are limited to Fourier bases,
then the computation of the coefficients must by necessity rely on
integration. \textquotedblleft Computers can't integrate!\textquotedblright\
Hmmm! Well, not directly. The problem must first be discretized. And there is
need for a more direct and algorithmic approach. Hence the wavelet algorithm!
This is a trend that has already found its way into textbooks on numerical
analysis (see, e.g., \cite{Coh03,BCD03}), and engineering (see, e.g.,
\cite{Stra00,StNg96}). In any case, algorithms
are central in mathematics even if you do not concern yourself with computers.
And it is the engineering connections that inspired the most successful
algorithms
in our subject. Our thoughts on the use of these trends in teaching have been
tried out by the author, and they are now in \cite{BrJo02b}, and more recently
in \cite{Jor06a}.\medskip

\section{\label{Section2}Case studies}

\subsection{\label{S2.1}Multiresolutions}

While finite or infinite families of nested subspaces
are ubiquitous in mathematics, and have been popular in Hilbert-space
theory for generations (at least since the 1930s{}), this idea was revived in
a different guise in 1986 by St\'{e}phane Mallat, then an engineering graduate
student; see \cite{Mal89}. In its adaptation to wavelets,
the idea is now referred to as the multiresolution method.

What made the idea especially popular in the wavelet
community was that it offered a skeleton on which various discrete algorithms
in applied mathematics could be attached and turned into wavelet constructions
in harmonic analysis. In fact what we now call multiresolutions
have come to signify a crucial link between the world of discrete wavelet
algorithms, which are popular in computational mathematics and in engineering
(signal/image processing, data mining,
etc.)\ on the one side, and on the other side continuous wavelet bases
in function spaces, especially in $L^{2}(\mathbb{R}^{d})$. Further, the
multiresolution
idea closely mimics how fractals are analyzed with the use of finite function systems.

But in mathematics, or more precisely in operator theory,
the underlying idea dates back to work of John von Neumann, Norbert Wiener,
and Herman Wold, where nested and closed subspaces
in Hilbert space were used extensively in an axiomatic approach to stationary processes,
especially for time series. Wold proved that any (stationary) time series can
be decomposed into two different parts: The first (deterministic) part can be
exactly described by a linear combination of its own past, while the second
part is the opposite extreme; it is \emph{unitary},
in the language of von Neumann.

John von Neumann's version of the same theorem is a pillar in operator theory.
It states that every isometry in a Hilbert space
$H$ is the unique sum of a shift isometry
and a unitary operator, i.e., the initial Hilbert space
$H$ splits canonically as an orthogonal sum of two subspaces $H_{s}$ and
$H_{u}$ in $H$, one which carries the shift operator, and the other $H_{u}$
the unitary part. The shift isometry is defined from a nested scale of closed spaces
$V_{n}$, such that the intersection of these spaces is $H_{u}$.

However, St\'{e}phane Mallat was motivated instead by the notion of scales of
resolutions
in the sense of optics. This in turn is based on a certain \textquotedblleft
artificial-intelligence\textquotedblright\ approach to vision and optics,
developed earlier by David Marr at MIT \cite{Mar82}, an approach which
imitates the mechanism of vision in the human eye.

The connection from these developments in the 1980s{} back to von Neumann is
this: Each of the closed subspaces $V_{n}$ corresponds to a level of resolution
in such a way that a larger subspace represents a finer resolution. Resolutions
are relative, not absolute! In this view, the relative complement of the
smaller (or coarser) subspace in larger space then represents the visual
detail which is added in passing from a blurred image to a finer one, i.e., to
a finer visual resolution.

This view became an instant hit in the wavelet
community, as it offered a repository for the fundamental father and the
mother functions, also called the scaling function $\varphi$, and the wavelet function
$\psi$. Via a system of translation and scaling operators,
these functions then generate nested subspaces,
and we recover the scaling identities which initialize the appropriate algorithms.

What results is now called the family of pyramid algorithms in wavelet analysis.
The approach itself is called the multiresolution approach (MRA) to wavelets.
And in the meantime various generalizations (GMRAs) have emerged.

In all of this, there was a second \textquotedblleft
accident\textquotedblright\ at play: As it turned out, pyramid algorithms
in wavelet analysis now lend themselves via multiresolutions,
or nested scales of closed subspaces, to an analysis based on frequency bands.
Here we refer to bands of frequencies as they have already been used for a
long time in signal processing.

Even though J. von Neumann and H. Wold had been using nested or scaled
families of closed subspaces
in representing past and future for time series, in 1989 S. Mallat, an
engineering graduate student at the time, found that this same idea applies
successfully to the representation of visual resolutions \cite{Mal89}.
And even more importantly, it offers a variety of powerful algorithms
for processing of digital images.

Now parallel to all of this, pioneers in probability theory had in fact
developed versions of the same refinement analysis. For example, in the theory
of martingales, consistency
relations may naturally be reformulated in the language of nested subspaces
in Hilbert space.

One reason for the success in varied disciplines of the same geometric idea is
perhaps that it is closely modeled on how we historically have represented
numbers in the positional number system;
see, e.g., \cite{Knu81}. Analogies to the Euclidean algorithm
seem especially compelling; see, e.g., \cite{SzFo70}.

\subsection{\label{S2.2}Fractals}

Intuitively, think of a fractal as reflecting similarity of scales such as is
seen in fern-like images that look \textquotedblleft roughly\textquotedblright%
\ the same at small and at large scales. While there may not be agreement
about a rigorous mathematical definition, Mandelbrot originally defined
\emph{fractals} as sets whose Hausdorff--Besicovich dimension
exceeded their topological dimension, but later accepted all self-similar,
self-affine, or quasi-self-similar sets as fractals. Moreover, even more generally, the
self-similarity could refer alternately to space, and to time. And further
versatility was added, in that flexibility is allowed into the definition of
\textquotedblleft similar.\textquotedblright\ 

In the book \cite{Jor06a}, our focus is more narrowly on the self-affine
variant (where computations are relatively simple); but in addition, we
encounter the fractal concept in other contexts, e.g., for \emph{measures},
and for probability processes (such as fractal Brownian motion).
Further, we have stressed examples more than the general theory.

The fractal concept for measures is especially agreeable in the self-affine
case, since affine maps
act naturally on measures. So for each fractal dimension $s$, there is a
corresponding $s$-fractal probability measure
$\mu=\mu_{s}$ which is the unique solution to a natural fixed-point equation,
one which depends on $s$. Moreover, we may then recover this way the spatial
fractal set itself as the support of this measure $\mu$.

As for $s$-fractal Brownian motion
(f\kern1.5ptBm), the \textquotedblleft$s$-fractal\textquotedblright\ feature
there refers to how the position $X_{t}$ at time $t$ of the
f\kern1.5ptBm-process transforms under scaling
of $t$: If time $t$ scales by $c$, then the respective distributions
before and after scaling
are related by the power-law $c^{s}$. Specifically, for all $t$, the finite
distributions calculated for $X_{ct}$ and for $c^{s}X_{t}$ coincide; see, e.g.,
\cite{JMR01}.

\subsection{\label{S2.3}Data mining}

The problem of how to handle and make use of large volumes of data is a
corollary of the digital revolution. As a result, the subject of data mining
itself changes rapidly. Digitized information (data) is now easy to capture
automatically and to store electronically \cite{HTK05}. In science, in
commerce, and in industry, data represents collected observations and
information: In business, there is data on markets, competitors, and customers
\cite{AgKu04a,AgKu04b}. In manufacturing, there is data for optimizing
production opportunities, and for improving processes \cite{Kus02,Kus05}. A
tremendous potential for data mining
exists in medicine \cite{KLD01,KDS05}, genetics \cite{ShKu04}, and energy
\cite{KuBu05}. But raw data is not always directly usable, as is evident by
inspection. A key to advances is our ability to \emph{extract information and
knowledge} from the data (hence \textquotedblleft data
mining\textquotedblright), and to understand the phenomena governing data sources.

Data mining
is now taught in a variety of forms in engineering departments, as well as in
statistics and computer science departments. One of the structures often
hidden in data sets is some degree of \emph{scale}. The goal is to detect and
identify one or more natural global and local scales in the data. Once this is
done, it is often possible to detect associated similarities of scale, much
like the familiar scale-similarity from multidimensional wavelets,
and from fractals. Indeed, various adaptations of wavelet-like algorithms
have been shown to be useful. These algorithms
themselves are useful in \emph{detecting} scale-similarities, and are
applicable to other types of pattern recognition. Hence, in this context,
generalized multiresolutions
offer another tool for discovering structures in large data sets, such as
those stored in the resources of the Internet. Because of the sheer volume of
data involved, a strictly manual analysis is out of the question. Instead,
sophisticated query processors based on statistical and mathematical
techniques are used in generating insights and extracting conclusions from
data sets. But even such an approach breaks down as the quantity of data grows
and the number of dimensions increases. Instead there is a new research area
(knowledge discovery in databases (KDD)) which develops various tools for
automated data analysis.

However, statistics is still at the heart of the problem of \emph{inference}
from the data. The widespread use of statistics, pattern recognition, and
machine-learning algorithms
is somewhat hindered in many areas by our ability to collect large volumes of
data. The next limitation in the subject arises when the data is too large to
fit in the main computer memory. As a result, we are faced with new issues,
e.g., quality of data, creative data analysis, and data transformation \cite{Kus01}.

Theory and hypothesis formation now becomes critical in our task of deriving
insights into underlying phenomena from the raw data. Various adaptations of
wavelet-like algorithms
have again proved useful in detecting scale-similarities, and in other types
of pattern recognition. Hence in this context wavelet
ideas offer another tool for discovering structures in vast data sets, such as
those in the resources of the Internet. And there are now a variety of such
effective Web mining tools in use.

Areas of data mining
include problems of representation, search complexity, and automated use of
prior knowledge to help in a data search. Thus we see the beginnings of a new
science for efficient inference from massive data sets.

\section{\label{Section3}Technology and the classroom}

While information communication technology (ICT) has advanced in leaps and
bounds in the past decade, it has not to the same degree impacted the
processes of learning and our classroom activities; at least not in striking
ways. Sure, you may say there is the Internet and there is Power Point, but
their direct effect in the class room has still been relatively modest. There
could be good reasons for that. The effect on what we teach and what students
learn outside the classroom has been much more striking.

As of yet, the direct benefits of ICT to classroom learning have not been
documented in convincing ways. But it is worth keeping in mind that in the
past few years, the impact of ICT has referred to both the subjects presented
in class, and to the way they are presented.

Examples:

\subsection{\label{S3.1}Technology impacts the form and manner of teaching}

The past few years have witnessed a substantial and direct use of the internet
as part of classroom presentations. Examples: Java scripts with moving frames
illustrating algorithms, direct projection of material in books and in
software, Powerpoint, and other such software tools.

\subsection{\label{S3.2}Technology impacting the substance and the subject
taught}

There is now more emphasis on the teaching of algorithms and approximation. We
include more numerical illustrations, more images, more interdisciplinary math
(combining ideas from math, from engineering, and from CS, --physics too!)
Other trends are two fold:

\begin{enumerate}
\item \label{S3.2(a)}New topics: e.g., discrete wavelet algorithms.

\item \label{S3.2(b)}Old topics in a new light: e.g., Signal and image processing.
\end{enumerate}

\subsection{\label{S3.3}Difficulties in adaptation of various technologies to
the classroom}

As ICT is adapted to teaching there have been various difficulties, some
intrinsic to teaching (\ref{S3.3(a)}) and some to the infrastructure and
organization (\ref{S3.3(b)}):

\begin{enumerate}
\item \label{S3.3(a)}Funding shortages, poor understanding of what works and
what does not, rigid policies enforced by bureaucrats who have limited
knowledge about the subjects affected, a lack of common sense in policies
regarding the various implementations. Misuses and overuse of ICT have on
occasion had the unintended effect of putting the students to sleep!

It has worked better when teachers have paid close attention to how students
in fact use these tools themselves.

\item \label{S3.3(b)}Organizational development issues have not always been
addressed effectively in implementations of ICT in education. When asked, I
tend to warn against models that are too rigid. In fact rigid and centralized
policies have often backfired. In my experience, ICT have worked best when
teachers listen and pay attention to how students use these tools themselves.
Teachers should inspire, and not \textquotedblleft
force-feed\textquotedblright\ students! I often notice that the best teachers
are also the best listeners!
\end{enumerate}

\noindent\textbf{Acknowledgments.} The work reported and summarized here was
supported by grants from the U.S. National Science Foundation (DMS-0139473
(FRG), DMS-0457581). In addition the author gratefully acknowledges many
discussions over time with numerous colleagues and students, too many to list.
However, special mention should be made of two former Ph.D. students, Dorin
Dutkay and Myung-Sin Song, and mathematics colleagues Professors Akram
Aldroubi, Tom Branson, David Larson, Wayne Lawton, Paul Muhly, and Judy
Packer, as well as engineering colleagues Professors Andrew Kusiak and Albert Ratner.

\ifx\undefined\bysame
\newcommand{\bysame}{\leavevmode\hbox to3em{\hrulefill}\,}
\fi

\end{document}